\documentclass{article}
\usepackage{amsmath}
\usepackage{amssymb}
\usepackage{amsthm}
\usepackage{bbm}
\usepackage{graphicx}
\usepackage{psfrag}
\usepackage{verbatim}
\usepackage{url}
\usepackage[usenames,dvipsnames]{xcolor}
\usepackage{todonotes}
\usepackage{hyperref}
\usepackage{authblk}

\hypersetup{
	colorlinks=true,
	linkcolor=blue!80!black,
	urlcolor=blue,
}

\newcommand{\bigO}{\mathcal{O}}
\newcommand{\Expect}{\mathbb{E}}

\newcommand{\Pinter}{P_{\text{inter}}}
\newcommand{\PKem}{P^{\text{Kem}}}
\newcommand{\KemSet}{\text{Kem}}
\newcommand{\colSet}{\mathcal C}
\newcommand{\taumix}{\tau_{\text{mix}}}
\newcommand{\taurel}{\tau_{\text{rel}}}
\newcommand{\distTV}[2]{||#1 - #2||_{\text{TV}}}
\newcommand{\unit}{\mathbbm{1}}

\DeclareMathOperator{\tw}{tw}
\DeclareMathOperator{\pw}{pw}

\DeclareMathOperator{\pwstart}{start}
\DeclareMathOperator{\pwend}{end}

\newtheorem{theorem}{Theorem}
\newtheorem{proposition}{Proposition}
\newtheorem{lemma}{Lemma}
\newtheorem{claim}{Claim}
\newtheorem{corollary}{Corollary}
\newtheorem{remark}{Remark}
\newtheorem{conjecture}{Conjecture}

\title{Glauber dynamics for colourings of chordal graphs and graphs of bounded treewidth}
\author{Marc Heinrich\thanks{Work supported by EPSRC grants EP/S016562/1, “Sampling in hereditary classes”.}}
\affil{University of Leeds, UK}

\begin{document}
\maketitle

\begin{abstract}
	The Glauber dynamics on the colourings of a graph is a random process which consists in recolouring at each step a random vertex of a graph with a new colour chosen uniformly at random among the colours not already present in its neighbourhood. It is known that when the total number of colours available is at least $\Delta +2$, where $\Delta$ is the maximum degree of the graph, this process converges to a uniform distribution on the set of all the colourings. Moreover, a well known conjecture is that the time it takes for the convergence to happen, called the mixing time, is polynomial in the size of the graph. Many weaker variants of this conjecture have been studied in the literature by allowing either more colours, or restricting the graphs to particular classes, or both. This paper follows this line of research by studying the mixing time of the Glauber dynamics on chordal graphs, as well as graphs of bounded treewidth. We show that the mixing time is polynomial in the size of the graph in the two following cases:
	\begin{itemize}
		\item on graphs with bounded treewidth, and at least $\Delta +2$ colours, 
		\item on chordal graphs if the number of colours is at least $(1+\varepsilon) (\Delta +1)$, for any fixed constant $\varepsilon$.
	\end{itemize}
\end{abstract}

\section{Introduction}
The Glauber dynamics is a Markov Chain on spin systems of a graph. We are interested here in the particular case of colouring, for which the dynamics can be described as follows. Starting from an initial arbitrary colouring of a given graph $G$, at each step of the chain a vertex $v$ of $G$ and a new colour $c$ are chosen uniformly at random. If this new colour does not already appear in the neighbourhood of $v$, then $v$ is recoloured with the colour $c$, otherwise the colouring is not changed. This very simple process has been widely studied both in the statistical physics and computer science community. 

This process occurs naturally as an extremal case of a model for particle interaction, called the Potts model where particles represented by the vertices of a graph can choose different states with probabilities depending on the state of their neighbours. The Glauber dynamics on the colourings of a graph corresponds to the limit case when the temperature goes to zero of the antiferromagnetic version of the problem (when adjacent particles attempt to choose different states).

This random process is also interesting from other points of view. When it is ergodic, i.e., when the distribution converges to a stationary (uniform in our case) distribution, independently of its initial state, it makes for a very natural and extremely simple candidate for sampling random uniform colourings of a graph via Monte Carlo Markov Chain algorithms (MCMC). Moreover, it is well known (see \cite{JVV86} for example) that for many problems, sampling a random uniform solution, and counting approximately the number of solutions are polynomially equivalent problems. In other words, any polynomial time algorithm for one can be turned into a polynomial time algorithm for the other.

Unfortunately, if we want to build an efficient sampler, ergodicity is not sufficient but we need a stronger property called rapid mixing. This property states in informal terms that the time it takes for the process to reach a distribution sufficiently close to the stationary distribution is polynomial in the input (here the size of the graph). Hence, finding sufficient conditions for rapid mixing to hold is a key question to decide when Glauber dynamics is sufficient to obtain a random sampler, and when more complicated methods are necessary. This question is at the center of this paper's research and motivates our study.

\paragraph{Known results}
The Glauber dynamics on the colourings of graphs is a random process that has been widely studied in the literature. An important part of the research on this Markov Chain has been driven by the following conjecture:
\begin{conjecture}
	For any graph $G$, the Glauber dynamics on the $k$-colourings of a graph $G$ for $k \geq \Delta +2$ has polynomial mixing time.
\end{conjecture}
Despite being several decades old, this conjecture remains largely open, and much research has been directed at weaker versions of this conjecture where the graph is restricted to special classes of graphs and/or we allow for a larger number of colours. Remark that the condition $k \geq \Delta+2$ is necessary for some graphs, since the chain is not ergodic with $\Delta+1$ colours on a clique for example. One of the early results on the problem is a proof by Jerrum~\cite{Jer95} that the result holds if the constraint on the number of colours is $k > 2 \Delta$ instead of $k \geq \Delta +2$. This was later improved by Vigoda~\cite{Vig00} to $k \geq \frac {11}{6} \Delta$. This remained the best known result for general graphs until recently where this was improved~\cite{CDMPP19} to $k \geq (\frac {11}{6} - \varepsilon) \Delta$ for some $\varepsilon \approx 10^{-5}$.

On the other hand, many results have focused on improving the bound on the number of colours by adding additional constraints on the graphs. The case of graphs with large girth in particular has attracted a lot of attention~\cite{DF03, HV03, DFHV04, Mol04, HV06}. Other classes of graphs have been considered in the literature, such as planar graphs~\cite{HVV15}, graphs of logarithmic pathwidth~\cite{Var18} or random graphs~\cite{MS10, EHSV18}. More specific graphs such as grids~\cite{MAMB04, GMP04, GJMP06, Jal12}, regular trees~\cite{GJK10}, or hyperbolic tilings~\cite{BKMP05} have been considered in the literature due to their frequent occurrence in the statistical physics community. The case of trees (not necessarily balanced ones) is of specific interest since it is one of the most well understood cases. Indeed, it was shown in~\cite{LMP09} that the mixing time of the Glauber dynamics on trees is $n^{\Theta(\frac{\Delta}{k \log \Delta})}$ for any number of colours $k \geq 3$. This implies a polynomial mixing time when either the number of colours is at least $C \frac {\Delta}{\log \Delta}$ for some constant $\Delta$, or when the maximum degree of the tree is bounded by a constant. Since trees are a well understood class of graph from the point of view of the mixing time of the Glauber dynamics, it is natural to consider larger classes of graphs with a tree-like structure. This is the line of research that we pursue here by considering graphs of bounded treewidth and chordal graphs.

\paragraph{Our contribution} The contributions of this paper are twofold.

The first result concerns the mixing time of the Glauber dynamics on graphs of bounded treewidth. We show that the Glauber dynamics mixes in time $n^{\bigO(\tw(G)^3)}$, where $\tw(G)$ is the treewidth of the graph if the number of colours is at least $\Delta +2$. This result holds in the more general setting of list colouring, where each vertex $v$ is given a list of colours of size at least $\deg(v) +2$ and can only be coloured with a colour from the list. Our result can be compared with Proposition~1.1 from~\cite{BKMP05}, which proves that the Glauber dynamics has mixing time at most $(\Delta+1) n (k -1)^{\xi(G) +1}$, where $\xi(G)$ is the cut-width of the graph. Since the cut-width satisfies the inequality: $\xi(G) \leq \Delta \tw(G) \log (n)$, their result implies a polynomial mixing time when both the treewidth and the maximum degree of the graph are constant. In our result, the exponent depends only on the treewidth of the graph. In particular our result proves polynomial mixing even if there are vertices with large (i.e., unbounded) degree. The dependency of the exponent on the treewidth is slightly worse than the result of~\cite{BKMP05} however. Our result is also an improvement on the result of~\cite{Var18} which proves polynomial mixing on graphs of bounded treewidth when the number of colours is $(1+\varepsilon)(\Delta +1)$ for a fixed constant $\varepsilon$.
	
The second results proves that the Glauber dynamics on chordal graphs mixes in polynomial time when the number of colours is at least $(1 + \varepsilon)(\Delta +1)$ for some fixed constant $\varepsilon$. This second result proceeds by first comparing the Glauber dynamics with an other Markov chain, that we call Kempe dynamics, in which colourings are changed by flipping the colours in a Kempe chain. We prove that the two chains have similar mixing times, up to polynomial factors, and then bound the mixing time of the Kempe dynamics in a second step. Note that this second step works as soon as the number of colours is at least $\omega +1$, and only the first step (i.e., the comparison between the two chains) requires $(1+\varepsilon) (\Delta +1)$ colours.

\paragraph{Organisation} The paper is organized as follows. In Section~\ref{sec+prelim} we define some notations on graph theory and Markov chains that will be used throughout the paper. We also recall in Section~\ref{sec+mixing} a small number of known techniques that will be used in our proofs. In Section~\ref{sec+tw} we prove the polynomial upper bound on the mixing time for graph of bounded treewidth. We start in Section~\ref{sec+list} by proving some inequalities on the number of list colourings of a graph, before moving to the main argument of the proof in Section~\ref{sec+tech+tw}. Finally, in Section~\ref{sec+chordal} we prove that the mixing time of Glauber dynamics is polynomial for chordal graphs when the number of colours is at least $(1+\varepsilon)(\Delta +1)$. The proof of this result proceeds in two steps. A first step, which is described in Section~\ref{sec+compare} compares the mixing time of the Glauber dynamics with the Kempe dynamics which allows for a more general type of recolouring. A second step described in Section~\ref{sec+Kempe}, bounds the mixing time for the Kempe dynamics using a coupling argument.

\section{Preliminaries}
\label{sec+prelim}

\subsection{Graph Theory}

We start with some standard definitions from graph Theory. A graph $G = (V, E)$ is defined by a finite set of vertices $V$, and a set $E$ of edges (unordered pairs of vertices). Given two vertices $u$ and $v$ of a graph $G$, we denote by $uv$ the edge between these two vertices. Given a vertex $v$, we denote by $N_G(v) = \{w \in V, vw \in E \}$ the neighbours of $v$, and $\deg_G(v) = |N(v)|$ its degree. The maximum degree of the graph is denoted by $\Delta_G = \max_{v \in V} \deg(v)$. Whenever the graph $G$ is clear from the context, we will drop the subscript, and simply write $N(v)$, $\deg(v)$, and $\Delta$ for these respective notions. Given $S \subseteq V$ a subset of vertices, $G[S]$ denotes the subgraph induced by the vertices in $S$, i.e., the graph whose vertex set is $S$, and whose edges are the edges of $E$ with both endpoints in $S$. 

We are interested in a family of graphs called Chordal graphs which do not contain induced cycles of length $4$ or more. An other characterisation of these graphs is the following:
\begin{proposition}[\cite{FG65}]
	\label{prop+chordal}
	A graph $G$ is chordal if and only if there exists an ordering $v_1, \ldots, v_n$ of its vertices such that for all $i$, the neighbours of $v_i$ which appear before in the ordering induce a clique in $G$. Such ordering is called a perfect elimination ordering of the graph.
\end{proposition}

The treewidth of a graph is a measure of close to a tree the structure of a graph is. More precisely, a graph $G = (V, E)$ has a tree-decomposition of width $k$ if and only if there exists a tree $T$, such that each node $u$ of $T$ is labelled by a set of vertices $S_u \subseteq V$, called a bag, and satisfying the following constraints:
\begin{itemize}
	\item for every vertex $v$ of $G$, the set of nodes $u$ of $T$ such that $v \in S_u$ induces a subtree of $T$;
	\item for every edge $vw \in E$, there is a node $u$ of $T$ such that $\{v, w\} \subseteq S_u$;
	\item all the bags have size at most $k+1$.
\end{itemize}
The treewidth of a graph $G$, denoted $\tw(G)$ is the smallest $k$ such that $G$ admits a tree-decomposition of width $k$. Both Chordal graphs and graphs of bounded treewidth are hereditary: they are closed by taking induced subgraphs. Moreover, Chordal graphs can also be defined as the graphs which have a tree decomposition such that every bag induces a clique in the graph $G$. In this case the decomposition is called a clique tree of $G$. Similarly, a path decomposition is defined in the way as above, with the additional constraint that $T$ must be a path. The smallest $k$ such that $G$ has a path decomposition of width $k$ is called the pathwidth of the graph, and is denoted $\pw(G)$

Given a graph $G$ and a set of colours $\colSet$, a list assignment $L$ for $G$ is a function $L : V \rightarrow 2^\colSet$ which associates to each vertex $v$ a list of colours $L(v)$. Given a list assignment $L$, an $L$-colouring of $G$ is a function $\sigma : V \rightarrow \colSet$ which associates to each vertex $v$ a colour $\sigma(v)$ such that $\sigma(v) \in L(v)$ for all the vertices of the graph, and $\sigma(v) \neq \sigma(w)$ for every pair of adjacent vertices $(v, w)$. We denote by $\Omega_{G, L}$ the set of all $L$-colourings of $G$. Given an integer $k \geq 0$, the $k$-colourings of $G$ correspond to the special case where $|\colSet| = k$, and $L(v) = \colSet$ for all the vertices of the graph. In a similar way, we denote by $\Omega_{G, k}$ the set of $k$-colourings of~$G$. The smallest value $k$ such that $\Omega_{G,k}$ is not empty is called the chromatic number of the graph $G$, and is written $\chi(G)$. As before, to lighten the notations the subscripts will be dropped when the graph $G$ and the list assignment $L$ are clear from the context. Finally, $L$ is a $\deg + 2$-list assignment of $G$ if and only if $|L(v)| \geq \deg(v) +2$ for every vertex $v$ of the graph $G$.

If $\sigma$ is a partial $L$-colouring of a subset $S$ of $G$, then let $L^{S, \sigma}$ be the list assignment of $G - S$ such that any $L^{S, \sigma}$ colouring of $G - S$ can be extended to $S$ by assigning the colour $\sigma(v)$ to all the vertices $v \in S$. In other words $L^{S, \sigma}(v) = L(v) \setminus \sigma(N(v) \cap S)$. In the special case where $S$ is a single vertex $v$, and $\sigma(v) = c$, this list assignment will be denoted simply by $L^{v, c}$.

Finally, given a graph $G$, a $k$-colouring $\sigma$, and two colours $c_1$ and $c_2$, a Kempe chain is a connected subgraph of $G$ whose vertices are all coloured either $c_1$ or $c_2$ in $\sigma$ and which is maximal for inclusion. A Kempe exchange consists in swapping the two colours in a Kempe chain, and result in a new colouring of the graph $G$. Note that Kempe exchange are not valid for list colourings since it might break list constraints, but only works when all the lists are equal. Recolouring a single vertex is

\subsection{Markov Chains and Glauber Dynamics}

We will assume that the reader is familiar with basic terminology from Markov chain theory (see~\cite{LP17} for an introduction on the subject). Given a graph $G$, and an list assignment $L$ of $G$, the Glauber dynamics on the $L$-colourings of $G$ is a Markov Chain with state space $\Omega_{G, L}$, the set of all $L$-colourings of $G$, and whose transition are obtained by doing the following procedure at each step.
\begin{enumerate}
	\item Choose a vertex $v$ uniformly at random.
	\item Choose a colour $c$ uniformly among all the colours in $L(v)$.
	\item Recolour $v$ with the colour $c$ if it is not already used by one of the nieghbours of $v$.
\end{enumerate}
Written in more formal terms, the transition matrix $P_{G, L}$ can be defined formally as follows for all $\sigma, \eta \in \Omega_{G, L}$:
$$ P_{G, L}[\sigma \rightarrow \eta] = \begin{cases}
	\frac 1 {n|L(v)|} & \text{if $\sigma$ and $\eta$ only differ at the vertex $v$.} \\ 
	0 & \text{otherwise.} 
\end{cases}$$

We can easily see that this chain is symmetric, i.e., its transition matrix $P$ satisfies for every $\sigma, \eta \in \Omega_{G, L}$, $P[\sigma \rightarrow \eta] = P[\eta \rightarrow \sigma]$. Moreover, if $|L(v)| \geq \deg(v) + 2$ for all the vertices of the graph, then it well known that this chain is ergodic and converges to the uniform distribution $\pi_{G, L}$ on $\Omega_{G, L}$. We are interested in the speed of convergence of this chain to the stationary distribution.

Given two distribution $\mu$ and $\nu$, the total variation distance  $\distTV{\mu}{\nu}$ between the two distributions is defined as:

$$ \distTV{\mu}{\nu} = \frac 1 2 \sum_{x \in \Omega_{G, L}} |\mu(c) - \nu(x)| \;.$$

Given an ergodic Markov chain with transition matrix $P$ and state space $\Omega$, if $\pi$ is its stationary distribution and for any $x \in \Omega$ and $t \geq 0$, $\nu_x^t$ is its distribution after $t$ steps starting at the position $x$, then its mixing time is defined as:

$$ \taumix(P) = \min\left\{ t,\ max_{x \in \Omega} \distTV{\nu_x^t}{\pi}  \leq \frac 1 4 \right\}\;. $$

An other measure of the speed of convergence is the relaxation time $\taurel$ which is related to the eigenvalues of the transition matrix $P$. More specifically, we write $\lambda(P)$ the smallest (in absolute value) non-zero eigenvalue of $P- I$, where $I$ is the identity matrix, which is called the spectral gap of the matrix $P$, and define $\taurel(P) = \frac 1 \lambda$ its relaxation time. The relaxation time and the mixing time are related by the following inequality:

$$ \taumix(P) \leq \log\left(\frac 4 {\min_{x \in \Omega} \pi(x)}\right) \taurel(P)  \;.$$

In particular, in the case of the Glauber dynamics, we have for every graph $G$ the inequality $ \taumix(P_{G, L}) \leq \log(4|\Omega_{G,L}|) \taurel(P_{G, L})$. This means that for the purpose of finding polynomial upper bound on the mixing time, it is enough to bound the relaxation time of the chain instead.

\subsection{Mixing time upper bounds}
\label{sec+mixing}

This section describes three different known techniques for bounding the mixing time of a given Markov chain. The first two allows to compare a given chain, with one or several other which are hopefully more simple, while the last one provides directly an upper bound on the mixing time.

\paragraph{Projection/Restriction.} The first of these techniques was introduced in~\cite{JSPV04}, and we will refer to it as projection/restriction. Consider a reversible Markov chain with state space $\Omega$ and transition matrix $P$. On an informal level, this technique consists in partitioning the state space $\Omega$ into subsets $\Omega_i$ for $i \leq \ell$, and relates the relaxation time of the chain on the whole space to:
\begin{itemize}
	\item chains restricted to each subset $\Omega_i$, 
	\item a chain measuring the movement between the elements of the partition (the projection chain).
\end{itemize}
In more details, the restriction $P_i$ is the Markov Chain whose state space is $\Omega_i$, and which has the same transition as the original chain for states in this subset. The projection chain has state space $\{1, \ldots, \ell\}$ and has the following transition for any $i, j \leq \ell$:

$$ 
\bar P[i \rightarrow j] = \sum_{x \in \Omega_i} \sum_{y \in \Omega_j} \frac{\pi(x)}{\pi(\Omega_i)} P[x \rightarrow y] \;.
$$

Note that the stationary distribution of the projection chain is $\bar \pi(i) = \pi(\Omega_i)$, and the chain $\bar P$ is reversible.
We also denote by $\gamma$ the following quantity:
$$ \gamma = \max_{i \in I} \max_{x \in \Omega_i} \sum_{y \in \Omega \setminus \Omega_i} P(x, y). $$
This parameter measures the maximum probability of leaving the subset $\Omega_i$ in the original chain. If we denote by $\lambda_{\min}$ the minimum spectral gap for the restriction chains $P_i$ (assuming all these chains are ergodic), and $\bar \lambda$ the spectral gap for the projection chain, then the following result holds:

\begin{proposition}[\cite{JSPV04}]
	\label{prop+proj+restrict}
	We have the following upper bound on the spectral gap $\lambda$ of $P$:
	$$ \lambda \geq \min\left(\frac {\bar \lambda}{3}, \frac{\lambda_{\min} \bar \lambda}{3\gamma + \bar \lambda}\right).$$
\end{proposition}

\paragraph{Canonical Paths.} A second comparison technique, sometimes also called multi-commodity flow was introduced in~\cite{Sin92}, and allows to bound the relaxation time of a Markov chain by constructing a system of paths between states of the Markov chain satisfying certain properties. This method was later improved~\cite{DS93} to allow the comparison of an unknown chain $P$, with an hopefully simpler chain $P'$, provided that the transition in $P'$ can be simulated by the transitions in $P$ in a way that does not put too much congestion on any one transition of $P$. Note that both chains are assumed to be reversible for the method to work. More precisely, denote by $\pi$ and $\pi'$ the stationary distribution for $P$ and $P'$ respectively, and assume that for every transition $(\alpha, \beta)$ of $P'$, we can associate a set $\Gamma_{\alpha, \beta}$ of paths from $\alpha$ to $\beta$ using only transition of $P$. Each path $\gamma$ is given a weight $g(\gamma)$ such that $\sum_{\gamma \in \Gamma_{\alpha, \beta}} g(\gamma) \geq 1$. Then, given a transition $(\sigma, \eta)$ of $P$, the congestion at this transition is given by:

$$ \rho_{\sigma, \eta} = \frac 1 {\pi(\sigma) P[\sigma, \eta]} \sum_{\alpha, \beta \in \Omega} \sum_{\gamma \in \Gamma_{\alpha, \beta}} g(\gamma) |\gamma| \pi'(\alpha) P'[\alpha, \beta]\;, $$
where $|\gamma|$ is the length of the path $\gamma$. Then we have the following upper bound on the relaxation time of $P$:

\begin{proposition}[Theorem 2.3 from~\cite{DS93}]
	\label{prop+canon+path}
	The following holds:
	$$ \taurel(P) \leq \taurel(P') \cdot \max_{(\sigma, \eta) \in P} \rho_{\sigma, \eta} \;.$$
\end{proposition}

A very special case of the result above is when the two Markov chains $P$ and $P'$ have the same transition, but with different transition probabilities. In this case, the transitions of $P'$ can be simulated directly by the transitions of $P$ using paths of length $1$.
\begin{corollary}
	\label{cor+canonPath}
	Given two Markov chains with transition probabilities $P$ and $P'$ on a state space $\Omega$; if there is a constant $C$ such that for every transition $(x, y) \in P'$, we have $\frac{P[x \rightarrow y]}{P'[x \rightarrow y]} \geq \frac 1 C$, then $\taurel(P) \leq C \cdot \taurel(P')$.
\end{corollary}

\paragraph{Coupling} A coupling of a Markov chain is a random process $(X_t, Y_t)_{t \geq 0}$ such that both $(X_t)_{t\geq 0}$ and $(Y_t)_{t \geq 0}$ are Markov chains with the same transition matrix $P$. Hence a coupling is defined by the joint distribution between two instances of the same Markov chain (with different initial states). The coupling time $T$ is the random variable whose value $t$ is the first time that $X_t = Y_t$. A coupling can be used to bound the mixing time of a Markov chain using the result below, which follows immediately from Theorem~5.2 from ~\cite{LP17} together with Markov's inequality.

\begin{proposition}
	\label{prop+coupling}
	If $(X_t, Y_t)_{t \geq 0}$ is a coupling for a Markov chain on state space $\Omega$ with transition matrix $P$ and stationary distribution $\pi$, and $T$ is the coupling time, then the following inequality holds:
	$$ \taumix(P) \leq 4 \max_{x, y \in \Omega} \Expect[T| X_0 = x, Y_0 = y]\;. $$
\end{proposition}

\section{Graphs of bounded treewidth}
\label{sec+tw}

The goal of this section is to prove the first of our two results, concerning graphs of bounded treewidth. Our result holds for the more general case of $(\deg +2)$-list-colouring. The exact statement is the following:
\begin{theorem}
	\label{thm+treewidth}
	Given a graph $G$, and a $(\deg +2)$-list assignment $L$ of $G$, the Glauber dynamics on the $L$-colourings of $G$ has mixing time at most $n^{\bigO(\tw(G)^3)}$.
\end{theorem}

In the rest of this section, $\lambda(G, L)$ denotes the spectral gap for the Glauber dynamics on the $L$-colourings of $G$, and $\lambda(G)$ is the minimum gap over all possible $(\deg+2)$-list assignments of $G$. In a similar way, we denote $\taurel(G) = \frac 1 {\lambda(G)}$ the maximum relaxation time over all possible list assignments. Let us denote by $\alpha_{\chi(G)}$ the quantity $\alpha_{\chi(G)} = 24(\chi(G) - 1)$. Note that the exact value of $\alpha_{\chi(G)}$ is not very important in our proofs, all we need to remember is that $\alpha_{\chi(G)}$ is a constant when the treewidth of the graph $G$ is bounded. This quantity will appear in the constant of several of our lemmas. Theorem~\ref{thm+treewidth} will follow immediately from the following result:

\begin{lemma}
	\label{lem+delv}
	Let $G$ be a graph containing at least $2$ vertices, and $v$ a vertex of $G$. The following inequality holds:
	$$ \taurel(G) \leq C_{\chi(G)} \taurel(G - v)\;,$$
	where $C_{\chi(G)} = 2^{- 8 \alpha_{\chi(G)}^2}$
\end{lemma}

Before we prove the lemma, let us see how it implies Theorem~\ref{thm+treewidth}. 
\begin{proof}[Proof of Theorem~\ref{thm+treewidth}]
	Let $C_{\chi(G)}$ be the constant defined in Lemma~\ref{lem+delv}. We will show by induction on the size of the tree decomposition $T$ of $G$ that $\taurel(G) \leq n \cdot (C_{\chi(G)})^{(\tw(G) +1) (\log(|T|) +1)}$ where $n$ is the number of vertices in $G$. Given a node $u$ of $T$, we denote by $S_u$ the bag associated to $u$. In the base case where the tree decomposition $T$ contains a single node $u$, the result follows by applying Lemma~\ref{lem+delv} to each of the vertices in the bag of $S_u$ until a single vertex remains in which case the relaxation time is equal to $1$.
	
	For the general case, observe that there exists a node $u$ in the tree $T$ such that the connected components of $T-u$ have size at most $\frac {|T|}{2}$. Let us denote by $G_1, \ldots G_d$ the connected components of $G - S_u$, and $n_i$ be the size of $G_i$. By  definition, each $G_i$ admits a tree decomposition $T_i$ with $|T_i| \leq \frac {|T|} {2}$. By applying Lemma~\ref{lem+delv} on each of the vertices of $S_u$ successively, we obtain:
	$$ \taurel(G) \leq (C_{\chi(G)})^{(\tw(G) +1)} \taurel(G - S) \;.$$
	Since $G-S$ is disconnected, the Glauber dynamics on $G-S$ is a product chain, and we have additionally: $\taurel(G-S) \leq \max_{i \leq d} \frac{n}{n_i} \taurel(G_i)$ (see Corollary~12.12 from~\cite{LP17}). Hence by applying the induction hypothesis on each of the connected components $G_i$, the following result holds:
	\begin{align*}
		\taurel(G) &\leq (C_{\chi(G)})^{\tw(G) +1}  \max_{i \leq d}\left(\frac{n}{n_i} \taurel(G_i)\right) \\
		&\leq (C_{\chi(G)})^{\tw(G) +1} \max_{i \leq d}\left(n \cdot(C_{\chi(G)})^{(\tw(G) +1) (\log(|T_i|) +1)}\right) \\
		&\leq n \cdot (C_{\chi(G)})^{(\tw(G) +1) (\max_i\log(|T_i|) +2)} \\
		&\leq n \cdot (C_{\chi(G)})^{(\tw(G) +1) (\log(|T|) +1)}
	\end{align*}
	This concludes the induction step. Since we can assume that the tree decomposition $T$ contains at most $n$ nodes, it follows that:
	$$ \taurel(G) \leq n^{\bigO(\tw(G)^3)}\;. $$
	Since the relaxation time and the mixing time are equal up to polynomial factors, this proves the theorem.
\end{proof}

Hence, to complete the proof of Theorem~\ref{thm+treewidth}, we only need to prove Lemma~\ref{lem+delv}. Before proving this lemma, we will need a few technical results to bound the number of list colourings of a graph, which are proved the section below. These bounds will then be used to prove Lemma~\ref{lem+delv} in Section~\ref{sec+texch+tw}.

\subsection{Bounds on the number of list colourings}
\label{sec+list}

In this section we give two bounds on the number of list colourings of a graph. These bounds will rely on the following technical result:

\begin{lemma}
	\label{lem+tech+bound}
	Let $A$ and $k$ be two positive integers, with $k \geq A$, $\varepsilon$ be a positive real with $\varepsilon < 1$, and $x_1, \ldots x_k$ be reals such that $0 \leq x_i \leq \frac 1 - \varepsilon$. If $\prod_{i \leq k} (1 - x_i) \leq \varepsilon^{A}$, then $\sum_{i \in k} x_i \geq (1- \varepsilon) A $.
\end{lemma}
\begin{proof}
	Let $x_i$ be reals which minimize the quantity $\sum_{i \leq k} x_i$ under the conditions of the lemma. Without loss of generality, we can assume $x_1 \geq x_2 \geq \ldots \geq x_k$. If $k = A$, then the only way to achieve the condition $\prod_{i \leq k} (1 - x_i) \leq  \varepsilon^{A}$ is if $x_i = 1 - \varepsilon$ for all $i \leq k$. In this case, $\sum_{i \leq k} x_i = (1 - \varepsilon) k  = (1 - \varepsilon) A $, and the result holds.
	
	We will show that even when $k > A$, the case above is the worst case. Let us assume by contradiction that $k > A$, and $x_{A+1} > 0$. Then we must have $x_A < 1 - \varepsilon$, since if we had an equality then we could decrease the value of the sum by setting $x_A$ to zero, and the constraint would still be satisfied. Let $\eta > 0$ such that $\eta < 1 - \varepsilon - x_A$, and $\eta < x_{A+1}$. Let us consider the weights $x'_i$ with $x'_A = x_A + \eta - \eta^2$ and $x'_{A+1} = x_{A+1} - \eta$, and $x'_i = x_i$ for $i \neq A, A+1$. Then we sill have by construction $0 \leq x'_i \leq 1 - \varepsilon$. Moreover, the following holds:
	
	\begin{align*}
	\prod_{i \leq k} (1 - x'_i) &= (1 - x_A - \eta + \eta^2)(1 - x_{A+1} + \eta) \prod_{i \neq A, A+1} (1 - x_i) \\ 
	&= ((1 - x_A)(1-x_{A+1}) + \eta(x_{A+1} - x_A) - \eta^2(x_{A+1} - \eta))\prod_{i \neq A, A+1} (1 - x_i) \\ 
	&\leq \prod_{i \leq k} (1 - x_i) \\ 
	& \leq 2^{-A}
	\end{align*}
	Where the inequality on the third line comes from the fact that both $x_{A+1} - x_A$ and $\eta - x_{A+1}$ are non-positive.
	Since $\sum_{i \leq k} x'_i = - \eta^2 + \sum_{i \leq x} x_i$, this contradicts the assumption that the $x_i$ minimize the sum under the constraints of the Lemma. It follows that we must have $x_{A+1} = 0$, which implies that $x_i = 0$ for all $i \geq A+1$. Finally, since we must have $\prod_{i \leq k} (1 - x_i) \leq 2^{-A}$, we must have $x_i = 1 - \varepsilon $ for all $i \leq A$ for this condition to be satisfied. This implies immediately $\sum_{i \leq k} x_i \geq (A - \varepsilon) A$, which holds for all possible choice of $x_i$ satisfying the conditions, and proves the lemma.
\end{proof}

We can now state and proofs the two results we will use to bound the number of $L$-colouring of a graph.

\begin{lemma}
	\label{lem+bound+Gv}
	Let $G$ be a graph, and $L$ be a $\deg + 1$-list assignment of $G$. If $v$ is a vertex with $|L(v)| \geq \deg(v) +2$, then the following inequality holds
	$$ |\Omega_{G, L}| \geq \max\left(\frac {|L(v)|} {\alpha_{\chi(G)}}, 2\right) \cdot |\Omega_{G-v, L}|$$
\end{lemma}
Note that in the Lemma above, the fraction $\frac 1 {\chi(G)  - 1}$ could be replaced by $\frac{|S|}{\deg(v)}$, where $S$ is the largest independent set of $N(v)$. This more general statement will however not be used in our proofs. The condition $|L(v)| \geq 3$ is there to rule out the case where $G$ is composed of just two adjacent vertices with both the same list of length $2$, for which the ratio $\frac{|\Omega_{G, L}| }{|\Omega_{G-v, L}|}$ is equal to~$1$.

\begin{proof}
	First, observe that the inequality $|\Omega_{G, L}| \geq 2 |\Omega_{G-v, L}|$ follows immediately from the fact that for every colouring of $G-v$, there are at least two choices of colours for $v$ to extend the colouring to the whole graph. Hence we can focus on the other inequality.
	
	Before proving the result for arbitrary graphs, let us consider the more simple case where $G$ is a star. The general case will be derived from this special case.
	\begin{claim}
		The lemma holds if $G$ is a star, and $v$ is the central vertex of the star.
	\end{claim}
	\begin{proof}
	Let us denote by $u_1, \ldots, u_d$ the neighbours of $v$ which are the leaves of the star.	By assumption on the list assignment $L$, we know that $|L(v)| \geq d+1$, and for every $u_i$, $|L(u_i)| \geq 2$. Moreover, since $G -v$ contains no edges, its number of colourings is exactly:
	$$ |\Omega_{G - v, L}| = \prod_{i \leq d} |L(u_i)|\;.$$
	Additionally, since $G$ is a star the number of colourings of the whole graph is equal to:
	$$ |\Omega_{G, L}| = \sum_{c \in L(v)} \prod_{i \leq d} |L(u_i) \setminus c|\;. $$
	Hence, the ratio between the two quantities is equal to:
	\begin{align} 
	\frac{|\Omega_{G, L}|}{|\Omega_{G - v, L}|} = \sum_{c \in L(v)} \prod_{i \leq d} \left( 1 - \frac{\unit_{c \in L(u_i)} }{|L(u_i)|}\right) \;. \label{eq+ratio+colouring} 
	\end{align}
	Let us denote by $J$ the set of colours $c$ for which the product which appears in the equation above is small. More specifically, this set is defined by:
	$$ J = \left\{c \in L(v), \quad \prod_{i \leq d} \left( 1 - \frac{\unit_{c \in L(u_i)} }{|L(u_i)|}\right) \leq \frac 1 8 \right\} \;.$$
	Then the result for $G$ will follow immediately from the fact that there is a linear number of colours which are not in $J$. Indeed, we know by assumption that $|L(u_i)| \geq 2$, and for every colour $c$ in $J$, we know using Lemma~\ref{lem+tech+bound} with the parameters $\varepsilon = \frac 1 2$, $A = 3$, and $k = d$, that the following inequality holds:
	$$ \sum_{i \leq d} \frac{\unit_{c \in L(u_i)} }{|L(u_i)|} \geq \frac 3 2\;.$$
	If we sum the inequalities above for every colour $c$ in $J$, we obtain the following result:
	\begin{align*}
	\frac 3 2 |J| &\leq \sum_{c \in J} \sum_{i \leq d} \frac{\unit_{c \in L(u_i)} }{|L(u_i)|} \\ 
	&\leq \sum_{i \leq d} \sum_{c \in L(v)} \frac{\unit_{c \in L(u_i)} }{|L(u_i)|} \\ 
	&= d
	\end{align*}
	Hence, it follows that $|J| \leq \frac {2d} 3$, and by definition of $J$, and using equality~\ref{eq+ratio+colouring}, we have:
	$$\frac{|\Omega_{G, L}|}{|\Omega_{G - v, L}|} \geq \frac 1 8 (|L(v)| - |J|) \geq \frac 1 {24} |L(v)| \;.$$
	
	This concludes the proof of the claim when the graph $G$ is a star, and $v$ is the central vertex.
	\end{proof}

	Let us now consider the general case where $G$ is an arbitrary graph, and assume that $v$ is a vertex of $G$. We know that the neighbourhood of $v$ contains an independent set in $G$ of size at least $\frac{\deg(v)}{\chi(G) - 1}$, and let $S$ be this independent set plus the vertex $v$. By construction, the graph induced by $S$ is a star, and let us denote $G' = G - S$, and $V' = V \setminus S$. We can partition the set $\Omega_{G, L}$ depending on the value taken by the colouring on the subgraph $G'$. Using this remark, the following equality holds:
	
	\begin{align}
	|\Omega_{G, L}| = \sum_{\sigma \in \Omega_{G', L}} |\Omega_{G[S], L^{ V', \sigma}}| \;. \label{eq+OmegaGL}
	\end{align}
	Applying the same decomposition to $G-v$, gives the similar result:
	\begin{align}
	|\Omega_{G -v, L}| = \sum_{\sigma \in \Omega_{G', L}} |\Omega_{G[S -v], L^{ V', \sigma}}| \;. \label{eq+OmegaGvL}
	\end{align}
	Moreover, we know that for every vertex $u \in S \setminus v$, we have $|L^{V', \sigma}(u)| \geq \deg_{G[S]}(u) +1$, and $|L^{V', \sigma}(v)| \geq \deg_{G[S]}(v) +2$, and since $G[S]$ is a star, applying the inequality we just obtained for stars to each term of the sum in equality~(\ref{eq+OmegaGL}) gives us:
	\begin{align*}
	|\Omega_{G, L}| &= \sum_{\sigma \in \Omega_{G', L}} |\Omega_{G[S], L^{ V', \sigma}}| \geq \sum_{\sigma \in \Omega_{G', L}} \frac 1 {24} |L^{V', \sigma}(v)| \cdot |\Omega_{G[S - v], L^{V', \sigma}}| \\ 
	\end{align*}
	Moreover, the following inequality holds:
	\begin{align*}
	|L^{V', \sigma}(v)| &\geq |L(v)| - |N(v) \setminus S| \\
	&\geq |L(v)| - \deg(v) + \frac {\deg(v)}{\chi(G) - 1} \\ 
	&\geq |L(v)|\left(1 - \frac {\chi(G) - 2}{\chi(G) -1} \cdot \frac{\deg(v)}{|L(v)|}\right) \\ 
	&\geq \frac 1 {\chi(G) -1} |L(v)|
	\end{align*}
	Where the last inequality comes from the fact that $|L(v)| \geq \deg(v)$. Combining these last two inequalities, and using equality~(\ref{eq+OmegaGvL}) gives the desired result:
	\begin{align*}
	|\Omega_{G, L}| &\geq \frac {|L(v)|} {24 (\chi(G) -1)} \sum_{\sigma \in\Omega_{G', L}}  |\Omega_{G[S - v], L^{ V', \sigma}}| \\ 
	&= \frac {|L(v)|} {24 (\chi(G) -1)} |\Omega_{G-v, L}| \\ 
	\end{align*}
\end{proof}

The lemma above can be interpreted by saying that if we consider a random uniform colouring of $G-v$, then there are on average a linear fraction of the colour of $|L(v)|$ which can extend this colouring to the whole graph $G$. In particular, it implies if we consider a random $L$-colouring $\sigma$ of $G$, then for every colour $c$, the probability that $\sigma(v) = c$ cannot be too high.

\begin{corollary}
	\label{cor+bound+distrib}
	Let $G$ be a graph and $L$ a $(\deg +2)$-list assignment of $G$, then for every vertex $v$, and every colour $c \in L(v)$:
	$$ \frac{|\Omega_{G - v, L^{v, c}}|}{|\Omega_{G, L}|} \leq \min\left(\frac 1 2, \frac {\alpha_{\chi(G)}}{|L(v)|}\right) \;.$$
\end{corollary}
\begin{proof}
	The result follows immediately from the fact that $|\Omega_{G - v, L^{v, c}}| \leq |\Omega_{G - v, L}|$, and the result from Lemma~\ref{lem+bound+Gv}.
\end{proof}

If the lemma above states gives an idea of how the number of colourings change when the vertex $v$ is removed from the graph, the result below gives a similar result when the colour of the vertex $v$ is fixed to an arbitrary value instead. Note that given a property $X$, the notation $\unit_X$ is a value equal to $1$ if the property $X$ is true, and equal to zero otherwise.

\begin{lemma}
	\label{lem+countLB}
	Let $G$ be a graph, and $L$ a $(\deg 12)$-list assignment of $G$, then the following holds for every vertex $v$:
	\begin{align*}
	\frac {|\Omega_{G -v, L^{v, c}}| }{|\Omega_{G-v, L}|} \geq \prod_{w \in N(v)} \left(1 - \min\left(\frac{\alpha_{\chi(G)} \cdot \unit_{c \in L(w)}}{|L(w)|}, \frac 1 2\right) \right)
	\end{align*}
\end{lemma}
\begin{proof}
	Let $u_1, \ldots, u_d$ be the neighbours of $v$, with $d = \deg(v)$, and let us denote by $L = L_0, L_1, \ldots, L_d =L^{v, c}$ the list assignments such that $L_i$ is obtained from $L$ by removing the colour $c$ from the list of the vertices $u_1, \ldots, u_i$. Then the ratio that we wish to bound can be written as the following telescopic product:
	$$ \frac {|\Omega_{G -v, L^{v, c}}| }{|\Omega_{G-v, L}|} = \prod_{i < d} \frac{|\Omega_{G-v, L_{i+1}}|}{|\Omega_{G - v, L_i}| \;.}$$
	Hence, to obtain the result, we only need to bound each of the terms separately, which is equivalent to proving that the result of the Lemma holds when $v$ has a single neighbour $w$. Notice that in this case $L^{v, c}$ is the list assignment where the colour $c$ has be removed from $L(w)$. If $L(w)$ does not contain the colour $c$, then $L^{v, c} = L$ which means that the ratio is equal to $1$, and the result holds. On the other hand, if $c\in L(w)$, then the ratio $ \frac {|\Omega_{G -v, L^{v, c}}| }{|\Omega_{G-v, L}|}$ represents the fraction of colourings of $\Omega_{G-v, L}$ which do not assign the colour $c$ to the vertex $w$. Hence, we obtain the following result:

	\begin{align*}
	\frac {|\Omega_{G -v, L^{v, c}}| }{|\Omega_{G-v, L}|} &= 1 - \frac {|\Omega_{G -v - w, L^{w, c}}| }{|\Omega_{G-v, L}|} \\ 
	&\geq 1 - \frac{|\Omega_{G -v - w, L}|}{|\Omega_{G-v, L}|} \\
	&\geq 1  - \min\left(\frac {24(\chi(G) - 1)} {|L(w)|}, \frac 1 2 \right) 
	\end{align*}
	Where the first equality comes from the fact that $w$ is $v$'s unique neighbour, the second inequality immediately follows from the fact that $L^{w,c}(u) \subseteq L(u)$ for every vertex $u$. The last inequality is obtained by applying Lemma~\ref{lem+bound+Gv} to the graph $G -v$, which can be applied since $|L(w)| \geq \deg_{G}(w) +1 = \deg_{G - v}(w) +2$.
\end{proof}

\subsection{Proof of Lemma~\ref{lem+delv}}
\label{sec+texch+tw}

We now have all the tools that we need to prove Lemma~\ref{lem+delv} and complete the proof of the theorem. Throughout this section we fix a graph $G$, with $v$ a vertex of $G$, and $L$ a $(\deg+2)$-list assignment of $G$. Our goal is to apply the restriction/projection technique by partitioning the set of $L$-colourings depending of the value the colourings take on the vertex $v$. More precisely, we can partition $\Omega_{G,L}$ into the sets $(\Omega_c)_{c \in L(v)}$, where $\Omega_c$ are the colourings which assign the colour $c$ to the vertex $v$. Note that there is a natural bijection between $\Omega_c$  and $\Omega_{G - v, L^{v, c}}$, and we will abuse notations and identify the two sets in the rest of this section. For this partition, the restriction chain has the following transition:
\begin{itemize}
	\item Pick a vertex $u$ uniformly at random, and a colour $c$ uniformly in $L(u)$.
	\item If $u \neq v$ and $c$ is not used by the neighbours of $u$, then recolour $u$ with the colour $c$.
\end{itemize}
From this description, it is clear that any restriction chain has almost the same transitions as the Glauber dynamics on $G-v$ for the list assignment $L^{v, c}$ for some colour $c$, with the two following modifications:
\begin{itemize}
	\item All the transition probabilities are reduced by a factor $\frac {n-1}{n}$ due to the possibility of choosing the vertex $v$  which does not modify the current colouring.
	\item For all the vertices $u \in N(v)$, the transition probabilities for colourings which differ at the vertex $u$ are further reduced by a factor of $\frac{|L(u)| - 1}{|L(u)|}$ due to the fact that the colour $c$ fixed by $v$ is still permitted in the choices for the new colour of $u$, but will result again in not modifying the colouring.
\end{itemize}
Since $L^{v, c}$ is a $(\deg+2)$-list assignment for the graph $G-v$, using Corollary~\ref{cor+canonPath} the observation above leads to the following remark:
\begin{remark}
\label{rem+restrict}
The spectral gap for any of these restriction chain is at least $\frac {\lambda(G-v)} 4$.
\end{remark}

Let us now consider the projection chain $\bar P$. The state space of the projection chain is $L(v)$, and the following transition probabilities holds for every pair of colours $(c_1, c_2)$ with $c_1 \neq c_2$: 

$$
\bar P[c_1 \rightarrow c_2] = \frac 1 {n|L(v)|} \frac {|\Omega_{G - v, L^{v, c_1}} \cap \Omega_{G - v, L^{v, c_2}}|} {|\Omega_{G - v, L^{v, c_1}}|}
$$

The stationary distribution of this Markov chain is $\bar \pi$, defined for every colour $c\in L(v)$ by $\bar \pi(c) = \pi(\Omega^{v, c})$. Moreover, we denote by $\gamma$ the following quantity:
$$ \gamma = \max_{c \in L(v)} \max_{\sigma \in \Omega_{G, L^{v, c}}} \sum_{y \in \Omega_{G, L} \setminus \Omega_{G, L^{v, c}}} P_{G, L}(x, y) \leq \frac 1 n.$$

In order to prove Lemma~\ref{lem+delv}, it will be enough to show the following result:

\begin{lemma}
	\label{lem+proj}
	The projection chain has spectral gap $\bar \lambda \geq \frac C n$, where $C = 2^{-7\alpha_{\chi(G)}^2}$.
\end{lemma}

Indeed, using the Lemma above and Remark~\ref{rem+restrict}, we can apply the projection restriction technique to get obtain the desired result.

\begin{proof}[Proof of Lemma~\ref{lem+delv}]
	We apply the projection/restriction technique with the partition described above. By Proposition~\ref{prop+proj+restrict}, the spectral gap of the Glauber dynmaics on the whole graph satisfies the following inequality:
	$$ \lambda(G, L) \geq \min\left(\frac{\lambda_{\min}}{3}, \frac {\lambda_{\min} \bar \lambda}{3 \gamma + \bar \lambda}\right)\;, $$
	where $\lambda_{\min}$ is the minimum spectral gap of one of the restriction chains, and $\bar \lambda$ is the spectral gap of the projection chain. By Remark~\ref{rem+restrict}, we know that $\lambda_{\min} \geq \frac {\lambda(G-v)}{4}$, hence it follows that there is a constant $C$ such that:
	\begin{align*}
	\frac {\lambda_{\min} \bar \lambda}{3 \gamma + \bar \lambda} &\geq \frac {\lambda(G-v)} 4 \cdot \frac {\frac C n}{\frac 3 n + \frac C n} \\ 
	& \geq \lambda(G-v) \cdot \frac {C} {4(3 + C)}
	\end{align*}
	The result immediately follows.
\end{proof}

Hence to complete the proof it only remains to prove Lemma~\ref{lem+proj}. The proof proceeds in several steps. The first step consists in defining a pair of colour $(c_1, c_2)$ as a good pair if $\frac {\bar{P}[c_1 \rightarrow c_2]}{\bar \pi(c_2)} \geq \frac{K}{n}$ for a constant $K = \frac {2^{- 6 \alpha_{\chi(G)}}}{\alpha_{\chi(G)}}$. The reason behind the value of the constant $K$ will become apparent later on.

We can now define an intermediate dynamics $\Pinter$ which has the following transitions: at each step, if the current state is the colour $c$, we choose a new colour $c'$ according to the distribution $\bar \pi$, and set $c'$ as the new state if $(c, c')$ is a good pair. In other words, this is the trivial chain, where a new state is chosen directly according to the stationary distribution, but with the restriction that  transitions are only allowed between good pairs of colours. We then use standard comparison arguments to compare first the spectral gaps of $\bar P$ and $\Pinter$, and then $\Pinter$ with the trivial chain which picks at each step a new state according to the stationary distribution $\bar \pi$.

\begin{lemma}
	\label{lem+comparePPinter}
	If $\lambda_{\text{inter}}$ is the spectral gap of $\Pinter$, then the following holds:
	$$ \frac{\bar \lambda}{\lambda_{\text{inter}}} \geq \frac{K}{n}$$
\end{lemma}
\begin{proof}
	By definition of a good pair, for every transition $(c_1, c_2)$ of $\Pinter$, we have $\bar P[c_1 \rightarrow c_2] \geq \frac K n \Pinter[c_1 \rightarrow c_2]$. Hence, the result follows immediately by Corollary~\ref{cor+canonPath}.
\end{proof}

Finally, we can bound the mixing time of the intermediate dynamics by showing that the fraction of pairs of colours which are not good pairs is sufficiently small.
\begin{lemma}
	\label{lem+Pinter}
	The spectral gap of $\Pinter$ satisfies $\lambda_{\text{inter}} \geq \frac 1 6$.
\end{lemma}
\begin{proof}
	We will again use a canonical path argument to prove the result. However, before describing how to construct the canonical paths, we will need an other result on the number of good pairs. Namely, we want to show that most of the pairs of colours are good pairs. 
	\begin{claim}
		For any colour $c$, if $\mathcal C_c$ is the set of colours $c'$ such that $(c, c')$ is a good pair, then $ \bar \pi(\mathcal C_c) \geq \frac 2 3$.
	\end{claim}
	\begin{proof}
		We know that:
		 $$ \bar{P}[c \rightarrow c'] = \frac 1 {n |L(v)|} \frac {|\Omega_{G - v, L^{v, c}} \cap \Omega_{G - v, L^{v, c'}}|} {|\Omega_{G - v, L^{v, c}}|} \;.$$
		 Moreover, the set $\Omega_{G - v, L^{v, c}} \cap \Omega_{G - v, L^{v, c'}}$ can also be written as $\Omega_{G - v, L'}$, where $L'$ is the list assignment obtained from $L$ by removing both $c$ and $c'$ from the lists of vertices adjacent to $v$. We can apply Lemma~\ref{lem+countLB} on the graph $G$ with the list assignment $L^{v, c}$. Indeed, since $L$ is a $\deg +2$-list assignment of $G - v$, and $L^{v, c}$ is obtained by removing a single colour, it is a $\deg+1$-list assignment of $G - v$. This gives the following inequality:
		\begin{align}
			\label{eq+lbTransProba}
			\bar{P}[c \rightarrow c'] &\geq \frac 1 {n |L(v)|} \prod_{w \in N(v)} \left(1 - \min\left(\frac 1 2, \frac {\alpha_{\chi(G)} \cdot \unit_{c' \in L(w)}}{|L(w) \setminus c|}\right)\right)
		\end{align}
		
		Let $A$ be a constant equal to 6 $\alpha_{\chi(G)}^2$. Let us first assume that $v$ has degree at most $A$. Then it follows that $\bar{P}[c \rightarrow c'] \geq \frac {2^{-A}} {n |L(v)|}$. Moreover, we know by Corollary~\ref{cor+bound+distrib} that $\bar \pi(c) \leq \frac{\alpha_{\chi(G)}}{|L(v)|}$. 
		Hence, this implies that $\frac{\bar{P}[c \rightarrow c']}{\bar {\pi}(c)} \geq \frac {2^{-A}}{n \alpha_{\chi(G)}} = \frac K n$. This means that all the pairs are good pairs, and there is nothing to prove.
		
		If we assume now that $v$ has degree at least $A$, then let $J$ be the set of colours $c'$ such that $(c, c')$ is not a good pair. By definition, this means that $\frac{\bar{P}[c \rightarrow c']}{\bar \pi(c)} \leq \frac K n$, and since $\bar \pi(c) \leq \frac{\alpha_{\chi(G)}}{|L(v)|}$, this implies that $\bar{P}[c \rightarrow c'] \leq \frac K n \bar{\pi}(c) \leq \frac{2^{-A}}{n |L(v)|}$. Hence, using this inequality, together with the inequation~\ref{eq+lbTransProba}, it follows that for any colour $c' \in J$, we have:
		$$ \prod_{w \in N(v)} \left(1 - \min\left(\frac 1 2, \frac {\alpha_{\chi(G)} \cdot \unit_{c' \in L(w)}}{|L(w) \setminus c|}\right)\right) \leq 2^{-A} \;.$$
		Hence, by Lemma~\ref{lem+tech+bound}, this implies that for $c' \in J$ we have:
		$$ \sum_{w \in N(v)} \frac {\alpha_{\chi(G)} \cdot 1_{c' \in L(w)}}{|L(w) \setminus c|} \geq \sum_{w \in N(v)} \min\left(\frac 1 2, \frac {\alpha \cdot \unit_{c' \in L(w)}}{|L(w) \setminus c|}\right) \geq \frac A 2$$
		By summing the inequality above for all the colours $c' \in J$, we obtain:
		$$ \alpha_{\chi(G)} \deg(v) \geq \frac A 2 |J| \;,$$
		which can be rewritten using the definition of $A$ as $|J| \leq \frac {2 \alpha_{\chi(G)}} {A} \deg(v) = \frac {1} {3 \alpha_{\chi(G)}} \deg(v)$. Additionally, we know by Corollary~\ref{cor+bound+distrib} that $\bar \pi(c_1) \leq \frac{\alpha_{\chi(G)}}{|L(v)|}$, which implies that $\bar{\pi}(J) \leq |J| \cdot \frac {\alpha_{\chi(G)}}{|L(v)|} \leq \frac {\deg(v)}{3 \alpha_{\chi(G)}} \cdot \frac {\alpha_{\chi(G)}} {|L(v)|} \leq \frac 1 3$, which proves the claim. 
	\end{proof}

	By claim above, for any two fixed colours $c_1$ and $c_2$, if $X$ is a colour chosen according to the distribution $\bar \pi$, then we have:
	
	\begin{align}
		\bar \pi(\mathcal C_{c_1} \cap \mathcal C_{c_2}) \geq 1 - \bar \pi(L(v) \setminus \mathcal C_{c_1}) - \bar \pi(\mathcal C_{c_2})  \geq \frac 1 3 \;. \label{eq+proba+goodPair}
	\end{align} 
	
	We can now build the canonical paths. We will build for each $c_1, c_2$ a set $\Gamma_{c_1, c_2}$ of weighted paths. For each $c' \in \mathcal C_{c_1} \cap \mathcal C_{c_2}$, we construct the canonical path $\gamma_{c'} = c_1, c', c_2$, with weight $g(\gamma_{c'}) = 3 \bar \pi(c')$. Note that we have as required $\sum_{c' \in \mathcal C_{c_1} \cap \mathcal C_{c_2}} g(\gamma_{c'}) = 3 \bar \pi(\mathcal C_{c_1} \cap \mathcal C_{c_2}) \geq 1$.
	
	Let us now compute the congestion for these canonical paths. Let $(c, c')$ a good pair, then:
	
	\begin{align*}
		\rho_{c, c'} &= \frac 1 {\bar \pi(c) \Pinter(c, c')} \sum_{c_1 \neq c_2} \sum_{\gamma \in \Gamma_{c_1, c_2}} g(\gamma) \bar \pi(c_1) \bar \pi (c_2) |\gamma| \\ 
		&= \frac 3 {\bar \pi (c) \bar \pi (c')} \left(\sum_{c_2, \gamma_{c'} \in \Gamma_{c, c_2}} g(\gamma_{c'}) \bar \pi(c) \bar \pi(c_2) + \sum_{c_1, \gamma_c \in \Gamma_{c_1, c'}} g(\gamma_c) \bar \pi(c_1) \bar \pi(c') \right) \\ 
		&\leq \sum_{c_2}  3\bar \pi(c_2) + \sum_{c_1} 3 \bar \pi(c_1) \\
		&\leq 6
	\end{align*}
	
	Since this holds for any pair of good colours $(c, c')$, the result follows immediately by applying Proposition~\ref{prop+canon+path}.
\end{proof}

\begin{proof}[Proof of Lemma~\ref{lem+proj}]
	The proof follows immediately from Lemma~\ref{lem+comparePPinter} and Lemma~\ref{lem+Pinter}.
\end{proof}

\section{Chordal graphs}
\label{sec+chordal}

In this section we prove the result below.

\begin{theorem}
	\label{thm+chordal}
	For any $\varepsilon > 0$, and for any chordal graph $G$ the Glauber dynamics on the $(1+\varepsilon) (\Delta + 1)$ colourings of $G$ has mixing time at most $n^{\bigO(1 + \log(1 + \frac 1 \varepsilon))}$.
\end{theorem}

The proof is quite different from the one for graphs of bounded treewidth. Instead of trying to decompose the graph into smaller parts, we extend the dynamics by adding additional transitions, and then use coupling argument to conclude. In more details, the proof of the theorem proceeds in two steps. The first step consists in comparing the Glauber dynamics with a new dynamics where a colouring can be modified in one step by applying a random Kempe exchange. Hence, the first step consists in using the canonical path technique to compare the mixing time of the two chains. This first step does not rely on the structure of chordal graph directly, but only on the property that Kempe chains on these graphs have bounded treewidth (they are in fact trees), and the fact that we have $\varepsilon\Delta$ additional colours. In a second step, we bound the mixing time for the Kempe dynamics. This is done using a simple coupling argument which uses the structure of chordal graphs. Note that the additional colours are not necessary for the second step, and the bound on the mixing time for the Kempe dynamics holds as soon as $k \geq \omega(G)$.

In the rest of this section, we denote by $P_{G, k}$ the Glauber dynamics on the graph $G$ with $k$ colours, where $k \geq \omega(G)$. We also define the dynamics $P^{\text{Kem}}_{G, k}$ the Kempe dynamics defined for any two $k$-colourings $\alpha \neq \beta$ of $G$ by:
\begin{align*}
	P^{\text{Kem}}_{G, k}(\alpha, \beta) = \begin{cases}
		\frac 1 {nk} & \text{if the two colourings differ by a Kempe exchange.} \\
		0 & \text{otherwise.}
	\end{cases}
\end{align*}
Note that the process above can be described with the following steps:
\begin{enumerate}
	\item select a vertex $v$ and colour $c$ uniformly at random; 
	\item let $C$ be the Kempe exchange which recolours $v$ with the colour $c$;
	\item swap the two colours in the Kempe chain with probability $\frac 1 {|C|}$.
\end{enumerate}

\subsection{Comparing $P_{G, k}$ and $\PKem_{G, k}$}
\label{sec+compare}

Our first step consists in comparing the relaxation time of the two Markov chains $P_{G, k}$ and $\PKem_{G, k}$:

\begin{lemma}
	\label{lem+compareKempe}
	For every $\varepsilon > 0$  and for every graph $G$, and $k \geq \Delta +1 + \varepsilon \omega$, we have:
	$$ \frac {\taurel}{\taurel^{\text{Kem}}} = n^{\bigO(1 + \log(1+ \frac 1 \varepsilon) )}$$
\end{lemma}
\begin{proof}
	The proof will follow by using canonical paths to simulate the transitions from $\PKem$ using only single vertex recolourings. Before we can describe how to build the canonical paths, we will need more details on how we use the structure of the graph. Since the graph $G$ is chordal, there is a pathwidth decomposition of $G$ such that every bag in the decomposition is the union of at most $\log(n+1)$ cliques (there can be edges between the cliques). The existence of such decomposition follows from the following claim:
	
	\begin{claim}
		Any chordal graph $G$ admits a path decomposition such that every bag can be partitioned into at most $\log(n) +1$ subsets each inducing a clique in $G$. 
	\end{claim}
	\begin{proof}
		We prove by induction on the size of the clique tree that every chordal graph $G$ admits a path decomposition where each bag can be partitioned into at most $\log(|T|) +1$ cliques, where $T$ is a clique decomposition of $G$. If the clique tree of $G$ contains a single node, then $G$ is a clique, and this clique tree is also a path decomposition which satisfy the constraints of the claim. 
		
		Let us now assume that $T$, the clique tree of $G$, has size at least~$2$. We know that there exists a node $u$ of $T$ such that the connected components of $T - u$ have size at most $\frac{|T|}{2}$. Let $S_u$ be the bag associated with the node $u$. By definition we know that $G[S_u]$ is a clique. Let $G_1, \ldots, G_\ell$ be the connected components of $G - S_u$. By construction, we know that each $G_i$ admits a clique tree of size at most $\frac{|T|}{2}$. Using the induction hypothesis, we know that for each $G_i$, there is a path decomposition such that each bag in the decomposition is the union of at most $\log\left(\frac{|T|}{2}\right) +1 = \log(|T|)$ cliques. The concatenation of all these path decomposition is again a path decomposition for $G - S_u$, and finally, by adding the set $S_u$ to all of the bags, we obtain a path decomposition of $G$ in which every bag is the union of at most $\log(|T|) +1$ cliques, which ends the induction and proves the claim.
	\end{proof}

	Let us denote by $B_1, \ldots B_\ell$ the bags in a pathwidth decomposition satisfying the constraints of the claim, and write $K = \log n +1$, the number of cliques each of the bags can be partitioned into. Every vertex $v$ is given two indices $\pwstart(v)$ and $\pwend(v)$ such that $v$ appears only in the bags with indices in the interval $[\pwstart(v), \pwend(v)]$.
	
	We now describe how to construct the canonical paths in order to compare $P_{G, k}$ and $\PKem_{G, k}$. Let $\alpha$ and $\beta$ be two colourings which differ by a Kempe chain $C$, with $|C| = t$. We now describe the set of paths $\Gamma_{\alpha, \beta}$, which is obtained as follows. For each $i$ from $1$ to $\ell$, at step $i$ the bag $B_i$ is recoloured in the following way:
	\begin{enumerate}
		\item Every vertex $v$ of $C \cap B_i$ with $\pwstart(v) = i$ is recoloured with an arbitrary colour different from $\alpha(v)$ (there are at least $k - \Delta - 1$ possible choices for the new colour).
		\item Every vertex $v$ of $C \cap B_i$ with $\pwend(v) = i$ is recoloured with its target colour $\beta(v)$.
	\end{enumerate}
	
	The fact that this procedure defines a valid recolouring path follows from the following two properties:\\
	\begin{itemize}
		\item At the end of step $i$, every vertex $v$ with $\pwend(v) \leq i$ is coloured with $\beta(v)$, and every vertex $v$ with $\pwstart(v) > i$ is coloured $\alpha(v)$.
		\item The vertices $v \in C$, with $\pwstart(v) \leq i \leq \pwend(v)$ are not coloured with their colour $\alpha(v)$.
	\end{itemize}
	The first property ensures that the transformation goes from $\alpha$ to $\beta$. Moreover, the two properties together ensure that we are never creating any monochromatic edge. Indeed, when a vertex $v$ is recoloured with its target colour $\beta(v)$, then all of its neighbours $w$ in $C$ either: (1) satisfy $\pwend(w) < i$, and consequently are already coloured according to $\beta$, or (2) $w \in B_i$, and by the second property, these vertices $w$ are not coloured with $\alpha(w) = \beta(v)$.
	
	Because there are at least $k - \Delta - 1$ choices the first time we recolour a vertex $v \in C$, this defines a family of paths $\Gamma_{\alpha, \beta}$ of size at least $(k - \Delta - 1)^t$ (recall that $t$ is the length of the Kempe chain). Hence, we can assign a weight of $(k - \Delta - 1)^t$ to each of these paths. To conclude the proof, it remains to bound the congestion associated to these canonical paths. Let $\sigma$ and $\eta$ be two colourings which differ on a single vertex $v$.

	Let us assume that:
	\begin{itemize}
		\item We know the index $i$ of the bag $B_i$ which is currently being processed by recolouring the vertex $v$ (there are at most $\ell \leq 2n$ possibilities).
		\item We know the two colours $c_1$ and $c_2$ which are involved in the Kempe chain $(\alpha, \beta)$ that we are considering (there are at most $k^2$ possiblities).
		\item We know the numbers $i_1$ and $i_2$ of vertices of $C$ coloured $c_1$ and $c_2$ in $\alpha$ in the bag $B_i$.
		\item We know which ones these vertices are. There are at most ${K \choose {i_1}} {K \choose {i_2}}\omega^{i_1 +i_2}$ possible choices for these vertices, obtained by first choosing which cliques these vertices are in, then picking one vertex in each of the cliques.
	\end{itemize}
	Note that once these choices are made, the colourings $\alpha$ and $\beta$ are uniquely defined. Inded, $\alpha$ can be obtained from $\sigma$ by recolouring the vertices of $B_i$ according to the choices made, and then performing Kempe exchange in $G[\bigcup_{j > i} B_j]$ between colours $c_1$ and $c_2$ as necessary to remove monochromatic edges between the vertices in $B_i$ and the rest of the graph. The colouring $\beta$ can be recovered in a similar fashion.
	
	Moreover, the fraction of paths of $\Gamma_{\alpha, \beta}$ which go through the transition $(\sigma, \eta)$ is at most $(k - \Delta - 1)^{-i_1 -i_2}$, since they are the paths which assign the correct colours for the vertices of $B_i$. Taking all this into account, it follows that:
	
	\begin{align*}
		\rho_{\sigma, \eta} &= \frac 1 {\pi(\sigma)P(\sigma, \eta)} \sum_{(\alpha, \beta) \in \PKem} \sum_{\gamma \in \Gamma_{\alpha, \beta}} g(\gamma) \pi(\alpha) \PKem(\alpha, \beta) \\ 
		&\leq n k^2 \sum_{i_1 \leq K} \sum_{i_2 \leq K} {K \choose {i_1}} {K \choose {i_2}}\omega^{i_1 +i_2} (k - \Delta - 1)^{-i_1 -i_2}\\
		&\leq n k^2 \sum_{i_1 \leq K} \sum_{i_2 \leq K} {K \choose {i_1}} {K \choose {i_2}}\left( \frac \omega {k - \Delta - 1}\right)^ {i_1 +i_2} \\
		& \leq n k^2 \left(1 + \frac 1 \varepsilon \right)^K \\
		&\leq n^{\bigO(1 + \log(1 + \frac 1 \varepsilon))}
	\end{align*}
	Where the fourth line follows from the bound on the number of colours, and the last line holds because $K = \log n +1$. The proof of the Lemma then follows immediately by applying Proposition~\ref{prop+canon+path}.
\end{proof}

Hence, to complete the proof of the theorem, it only remain to show that the Kempe dynamics has polynomial mixing time, which is proved in the next subsection.

\subsection{Kempe dynamics}
\label{sec+Kempe}

The upper bound on the mixing time for the Kempe dynamics is obtained using a coupling argument. We prove the following result:

\begin{theorem}
	\label{thm+Kempe}
	The Kempe dynamics on the $k$-colourings of a chordal graph $G$, with $k \geq \omega$ has mixing time at most $\bigO(\omega n^2)$.
\end{theorem}

In the rest of this section, given a colouring $\sigma$ of a graph $G$, we will denote by $\KemSet(\sigma)$ the set of Kempe exchange which can be applied to colouring $\sigma$. In order to simplify notations, we will assume that the empty chain is added to $\KemSet(\sigma)$ the corresponding number of times in order to get $|\KemSet(\sigma)| = kn$. With this abuse of notation, the Kempe dynamics can be defined as choosing a Kempe exchange in $\KemSet(\sigma)$ uniformly at random, and swapping the two colours in this chain (if this is the empty chain, then nothing happens).

Let $G$ be a chordal graph. We will assume that the vertices $v_1, \ldots, v_n$ of $G$ are ordered in a perfect elimination ordering. In other words, for every index $i$, the neighbours of $v_i$ which have a smaller index induce a clique in the graph $G$. Moreover, we denote by $V_i = \{v_1, \ldots v_i\}$, with the convention that $V_0 = \emptyset$.

We will define a coupling $(X_t, Y_t)$ for the Kempe dynamics by simply pairing the Kempe exchange of $\KemSet(X_t)$ with the elements of $\KemSet(Y_t)$. The corresponding coupled dynamics consists in choosing a random pair, and modifying $X_t$ and $Y_t$ accordingly.

In order to define this pairing, let $i$ be the first index such that $X_t$ and $Y_t$ differ on $v_i$. We can partition the set $\KemSet(X_t)$ into two sets: 
\begin{itemize}
	\item $K^X_1$ denotes the Kempe exchanges which affect at least one of the vertices of $V_i$.
	\item $K^X_2$ denotes the Kempe exchanges which do not affect the vertices of $V_i$.
\end{itemize}
In a similar fashion, we define $K^Y_1$ and $K^Y_2$ the partition of $\KemSet(Y_t)$. Observe that there is a natural bijection $\phi$ between $K^X_1$ and $K^Y_1$, which allows us to pair the elements of this set in a natural way. In more details, the bijection can be described as follows:
\begin{itemize}
	\item If $C_1$ is a Kempe exchange of $X_t$ which recolour $u \in V_{i-1}$ with the colour $c$, then $\phi(C_1)$ is the Kempe exchange of $Y_t$ which also recolours $u$ with the colour $c$. Observe that because of the elimination ordering, both $C_1$ and $\phi(C_1)$ intersect $V_{i-1}$ in the same way. In other words, after applying the change $C_1$ to $X_t$, and $\phi(C_1)$ to $Y_t$, the two colourings still agree on $V_{i-1}$. Note however than one of $C_1$ or $\phi(C_1)$ might recolour $v_i$ but not the other.
	\item If $C_1$ is a Kempe exchange of $X_t$ which recolours $v_i$ with the colour $c$, but does not recolour any other vertex of $V_i$, then $\phi(C_1)$ is the Kempe exchange of $Y_t$ which recolours $v_i$ with the colour $c$. Observe that applying $C_1$ and $\phi(C_1)$ to $X_t$ and $Y_t$ respectively, result in two colourings which agree on $v_i$, and leave $V_{i-1}$ unchanged.
\end{itemize}

Hence, the elements of $K^X_1$ are paired with the elements of $K^Y_1$ according to the bijection described above. The elements of $K^X_2$ and $K^Y_2$ are paired in an arbitrary fashion. Using this pairing, the coupling $(X_t, Y_t)$ satisfies the following properties:

\begin{lemma}
	\label{lem+coupling}
	The coupling defined above satisfies:
	\begin{itemize}
		\item For all $j$, if $X_t|{V_j} = Y_t|_{V_j}$, then $X_{t+1}|{V_j} = Y_{t+1}|_{V_j}$.
		\item The probability that $X_{t+1}$ and $Y_{t+1}$ assign the same colour to $v_i$ is at least $\frac{k - \omega - 1}{kn}$.
	\end{itemize}
\end{lemma}
\begin{proof}
	The first point follows from the definition of the bijection $\phi$ which satisfies $C \cap V_{i-1} = \phi(C) \cap V_{i-1}$ for all the Kempe exchange $C \in K_1$. The second point is a consequence of the fact that if we choose a Kempe exchange $C$ which recolours $v_i$ with a colour which is not already taken by one of its neighbours which appears before it in the elimination ordering, then both $X_t$ and $Y_t$ will assign the same colour to $v_i$. Since there are at least $k - \omega - 1$ such colours, the result follows immediately.
\end{proof}

The lemma above implies the following result, which proves Theorem~\ref{thm+Kempe} using Proposition~\ref{prop+coupling}.
\begin{corollary}
	\label{cor+expectTime}
	Let $T$ be the first time $t$ that $X_t = Y_t$, then $\Expect(T) \leq \omega n^2$
\end{corollary}
\begin{proof}
	By the first point of Lemma~\ref{lem+coupling}, we know that the first index $i$ on which $X_t$ and $Y_t$ disagree can only increase with time. Let $T_i$ denotes the first time that $X_t$ and $Y_t$ agree on $V_i$, and remark that we can write: $T = \sum_{i < n} T_{i+1} - T_i$. Hence by linearity of the expectation, $\Expect(T) = \sum_{i < n} \Expect(T_{i+1} - T_i)$. Moreover, $\Expect(T_{i+1} - T_i)$ denotes the expected time for the two colourings to agree on an additional vertex, and by the second point of Lemma~\ref{lem+coupling}, it follows that $\Expect(T_{i+1} - T_i) \leq \frac{nk}{k - \omega - 1} \leq \omega n$. From this, it follows that $\Expect(T) \leq \omega n^2$. 
\end{proof}

\begin{proof}[Proof of Theorem~\ref{thm+chordal}]
	The proof follows immediately from Lemma~\ref{lem+compareKempe} and Theorem~\ref{thm+Kempe}.
\end{proof}

\section{Conclusion}
We have shown that the Glauber dynamics has polynomial mixing time for graphs of bounded treewidth and at least $\Delta+2$ colours, as well as for chordal graphs and at least $(1+\varepsilon)(\Delta+1)$ colours. One further direction of research could be to try to decrease further the number of colours in the case of graphs of bounded treewidth. It is known that the chain is ergodic if the number of colours is at least $\tw(G) +2$. However, it is known that this number of colours is not sufficient to ensure rapid mixing, since $\Omega(\frac {\Delta}{\log \Delta})$ colours are required in the case of a star for rapid mixing to occur. This counter example still leaves the possibility that rapid mixing could occur when the number of colours is of the order $\tw(G) \frac{\Delta}{\log \Delta}$. An other question arises when we compare our result with the work of~\cite{BKMP05}. In their paper, they also prove that polynomial mixing occurs for the ferromagnetic Ising model on graphs of bounded treewidth, but with an exponent which depends on the maximum degree of the graph. It could be interesting to investigate whether our method also extends to the Ising model, and whether it can be used to obtain a bound on the mixing time with an exponent which is independent from the maximum degree of the graph.

In the case of chordal graphs, the question whether rapid mixing occurs with $\Delta+2$ colours is still open. Since we already proved that Kempe dynamics are rapidly mixing on chordal graphs for any number of colours, it might be possible to prove this result by improving on the techniques which allow to compare the mixing time of the Glauber dynamics with that of the Kempe dynamics. Additionally, this method of comparing these two chains could possibly be useful for other classes of graphs where the Kempe chains have a particular structure. One important example is the case of linegraphs (or equivalently, edge colourings) for which all Kempe chains are either paths or cycles. Provided that the comparison still works on this class, it might be possible to prove rapid mixing results for a smaller number of colours by studying the Kempe dynamics instead of the Glauber dynamics.

\bibliographystyle{alpha}
\bibliography{bib}

\newcommand{\etalchar}[1]{$^{#1}$}
\begin{thebibliography}{AMMVB04}

\bibitem[AMMVB04]{MAMB04}
Dimitris Achlioptas, Mike Molloy, Cristopher Moore, and Frank Van~Bussel.
\newblock Sampling grid colorings with fewer colors.
\newblock In {\em LATIN 2004: Theoretical Informatics}, pages 80--89. Springer
  Berlin Heidelberg, 2004.

\bibitem[BKMP05]{BKMP05}
Noam Berger, Claire Kenyon, Elchanan Mossel, and Yuval Peres.
\newblock Glauber dynamics on trees and hyperbolic graphs.
\newblock {\em Probability Theory and Related Fields}, 131(3):311--340, 2005.

\bibitem[CDM{\etalchar{+}}19]{CDMPP19}
Sitan Chen, Michelle Delcourt, Ankur Moitra, Guillem Perarnau, and Luke Postle.
\newblock Improved bounds for randomly sampling colorings via linear
  programming.
\newblock In {\em Proceedings of the Thirtieth Annual ACM-SIAM Symposium on
  Discrete Algorithms}, pages 2216--2234. SIAM, 2019.

\bibitem[DF03]{DF03}
Martin Dyer and Alan Frieze.
\newblock Randomly coloring graphs with lower bounds on girth and maximum
  degree.
\newblock {\em Random Structures \& Algorithms}, 23(2):167--179, 2003.

\bibitem[DFHV04]{DFHV04}
Martin Dyer, Alan Frieze, Thomas~P Hayes, and Eric Vigoda.
\newblock Randomly coloring constant degree graphs.
\newblock In {\em 45th Annual IEEE Symposium on Foundations of Computer
  Science}, pages 582--589. IEEE, 2004.

\bibitem[DSC93]{DS93}
Persi Diaconis and Laurent Saloff-Coste.
\newblock Comparison theorems for reversible markov chains.
\newblock {\em The Annals of Applied Probability}, pages 696--730, 1993.

\bibitem[EH{\v{S}}V18]{EHSV18}
Charilaos Efthymiou, Thomas~P Hayes, Daniel {\v{S}}tefankovi{\v{c}}, and Eric
  Vigoda.
\newblock Sampling random colorings of sparse random graphs.
\newblock In {\em Proceedings of the Twenty-Ninth Annual ACM-SIAM Symposium on
  Discrete Algorithms}, pages 1759--1771. SIAM, 2018.

\bibitem[FG65]{FG65}
Delbert Fulkerson and Oliver Gross.
\newblock Incidence matrices and interval graphs.
\newblock {\em Pacific journal of mathematics}, 15(3):835--855, 1965.

\bibitem[GJK10]{GJK10}
Leslie~Ann Goldberg, Mark Jerrum, and Marek Karpinski.
\newblock The mixing time of glauber dynamics for coloring regular trees.
\newblock {\em Random Structures \& Algorithms}, 36(4):464--476, 2010.

\bibitem[GJMP06]{GJMP06}
Leslie~Ann Goldberg, Markus Jalsenius, Russell Martin, and Mike Paterson.
\newblock Improved mixing bounds for the anti-ferromagnetic potts model on z 2.
\newblock {\em LMS Journal of Computation and Mathematics}, 9:1--20, 2006.

\bibitem[GMP04]{GMP04}
Leslie~Ann Goldberg, Russell Martin, and Mike Paterson.
\newblock Strong spatial mixing for lattice graphs with fewer colours.
\newblock In {\em 45th Annual IEEE Symposium on Foundations of Computer
  Science}, pages 562--571. IEEE, 2004.

\bibitem[HV03]{HV03}
Thomas~P Hayes and Eric Vigoda.
\newblock A non-markovian coupling for randomly sampling colorings.
\newblock In {\em 44th Annual IEEE Symposium on Foundations of Computer
  Science, 2003. Proceedings.}, pages 618--627. IEEE, 2003.

\bibitem[HV06]{HV06}
Thomas~P Hayes and Eric Vigoda.
\newblock Coupling with the stationary distribution and improved sampling for
  colorings and independent sets.
\newblock {\em The Annals of Applied Probability}, 16(3):1297--1318, 2006.

\bibitem[HVV15]{HVV15}
Thomas~P Hayes, Juan~C Vera, and Eric Vigoda.
\newblock Randomly coloring planar graphs with fewer colors than the maximum
  degree.
\newblock {\em Random Structures \& Algorithms}, 47(4):731--759, 2015.

\bibitem[Jal12]{Jal12}
Markus Jalsenius.
\newblock Sampling colourings of the triangular lattice.
\newblock {\em Random Structures \& Algorithms}, 40(4):501--533, 2012.

\bibitem[Jer95]{Jer95}
Mark Jerrum.
\newblock A very simple algorithm for estimating the number of k-colorings of a
  low-degree graph.
\newblock {\em Random Structures \& Algorithms}, 7(2):157--165, 1995.

\bibitem[JSTV04]{JSPV04}
Mark Jerrum, Jung-Bae Son, Prasad Tetali, and Eric Vigoda.
\newblock Elementary bounds on poincar{\'e} and log-sobolev constants for
  decomposable markov chains.
\newblock {\em The Annals of Applied Probability}, 14(4):1741--1765, 2004.

\bibitem[JVV86]{JVV86}
Mark~R Jerrum, Leslie~G Valiant, and Vijay~V Vazirani.
\newblock Random generation of combinatorial structures from a uniform
  distribution.
\newblock {\em Theoretical computer science}, 43:169--188, 1986.

\bibitem[LMP09]{LMP09}
Brendan Lucier, Michael Molloy, and Yuval Peres.
\newblock The glauber dynamics for colourings of bounded degree trees.
\newblock In {\em Approximation, Randomization, and Combinatorial Optimization.
  Algorithms and Techniques}, pages 631--645. Springer, 2009.

\bibitem[LP17]{LP17}
David~A Levin and Yuval Peres.
\newblock {\em Markov chains and mixing times}, volume 107.
\newblock American Mathematical Soc., 2017.

\bibitem[Mol04]{Mol04}
Michael Molloy.
\newblock The glauber dynamics on colorings of a graph with high girth and
  maximum degree.
\newblock {\em SIAM Journal on Computing}, 33(3):721--737, 2004.

\bibitem[MS10]{MS10}
Elchanan Mossel and Allan Sly.
\newblock Gibbs rapidly samples colorings of g (n, d/n).
\newblock {\em Probability theory and related fields}, 148(1-2):37--69, 2010.

\bibitem[Sin92]{Sin92}
Alistair Sinclair.
\newblock Improved bounds for mixing rates of markov chains and multicommodity
  flow.
\newblock {\em Combinatorics, probability and Computing}, 1(4):351--370, 1992.

\bibitem[Var18]{Var18}
Shai Vardi.
\newblock Randomly coloring graphs of logarithmically bounded pathwidth.
\newblock In {\em Approximation, Randomization, and Combinatorial Optimization.
  Algorithms and Techniques (APPROX/RANDOM 2018)}. Schloss
  Dagstuhl-Leibniz-Zentrum fuer Informatik, 2018.

\bibitem[Vig00]{Vig00}
Eric Vigoda.
\newblock Improved bounds for sampling colorings.
\newblock {\em Journal of Mathematical Physics}, 41(3):1555--1569, 2000.

\end{thebibliography}
\end{document}